\documentclass[10pt,a4paper]{amsart}
\usepackage{graphics,multirow,array,amsmath,amssymb,epsf}

\newtheorem{thm}{Theorem}
\newtheorem{lem}[thm]{Lemma}

\newtheorem{cor}[thm]{Corollary}

\theoremstyle{definition}

\begin{document}
\title
{Knotted Hamiltonian cycles in linear embedding of $K_7$ into $\mathbb{R}^3$}
\author{Youngsik Huh}
\address{Department of Mathematics,
College of Natural Sciences, Hanyang University, Seoul 133-791,
Korea} \email{yshuh@hanyang.ac.kr}




\begin{abstract}
In 1983 Conway and Gordon proved that any embedding of the complete graph $K_7$ into $\mathbb{R}^3$ contains at least one nontrivial knot as its Hamiltonian cycle. After their work  knots (also links) are considered as intrinsic properties of abstract graphs, and numerous subsequent works have been continued until recently. In this paper we are interested in knotted Hamiltonian cycles in linear embedding of $K_7$. Concretely it is shown that any linear embedding of $K_7$ contains at most three figure-8 knots as its Hamiltonian cycles.
\end{abstract}

\maketitle


\section{Introduction}
In this paper we are interested in knots residing in spatial embeddings of graphs. Before describing our interests concretely we need to give necessary definitions.

\vspace{0.2cm}
A circle embedded into the Euclidean 3-space $\mathbb{R}^3$ is called a {\em knot}. Two knots $K$ and $K^{\prime}$ are said to be {\em ambient isotopic}, denoted by $K\sim K^{\prime}$, if there exists a continuous map $h:\mathbb{R}^3\times[0,1] \rightarrow \mathbb{R}^3$ such that the restriction of $h$ to each $t \in [0,1]$, $h_t:\mathbb{R}^3\times\{t\} \rightarrow \mathbb{R}^3$, is a homeomorphism, $h_0$ is the identity map and $h_1(K_1)=K_2$, to say roughly, $K_1$ can be deformed to $K_2$ without intersecting its strand. The ambient isotopy class of a knot $K$ is called the {\em knot type} of $K$. Especially if $K$ is ambient isotopic to another knot contained in a plane of $\mathbb{R}^3$, then we say that $K$ is {\em trivial}.

\begin{figure}[h]
{\epsfxsize=7cm \epsfbox{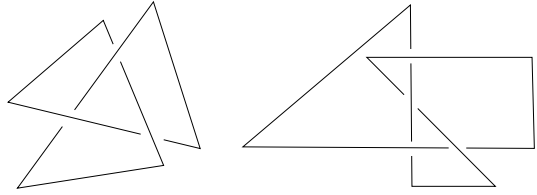}}
\caption{Polygonal presentations of $3_1$ and $4_1$ knots}
\label{fig1}
\end{figure}

A {\em polygonal knot} is a knot consisting of finitely many line segments. Figure \ref{fig1} shows polygonal presentations of two knot types $3_1$ and $4_1$ (These notations for knot types follow the knot tabulation in \cite{Rolfsen}. Usually $3_1$ and its mirror image are called {\em trefoil}, and $4_1$ {\em figure-8}).
Polygonal knots appear in many scientific contexts other than mathematics. Real knots such as polymers are modeled to be spatial  polygons for theoretical studies on their physical and chemical properties \cite{Randell-3, Sumners}. One quantity of polygonal knots interesting in the related research is polygonal index. For a knot type $\mathfrak{K}$, its {\em polygonal index} $p(\mathfrak{K})$ is defined to be the minimal number of line segments required to realize $\mathfrak{K}$ as a polygon. From the definition we easily know $p(\mbox{trivial knot})=3$. But generally it is not easy to determine $p(\mathfrak{K})$ for an arbitrary knot type $\mathfrak{K}$. This quantity was determined only for some specific knot types \cite{Calvo, FLS, J, Mc, Randell-2, CR}. Here we mention a result by Randell on small knots for later use.
\begin{thm} \label{stick-thm}
\cite{Randell-2}
$p(\mbox{trivial knot})=3$, $p(\mbox{trefoil})=6$ and $p(\mbox{figure-8})=7$. Furthermore, $p(\mathfrak{K})\geq 8$ for any other knot type $\mathfrak{K}$.
\end{thm}

Let $K_n$ be the complete graph with $n$ vertices. A spatial embedding of $K_n$, that is, an embedding of $K_n$ into $\mathbb{R}^3$, will be said to be {\em linear }, if each edge of the graph is mapped to a line segment. Note that a linearly embedded $K_n$ contains $\frac{n!}{2n}$ polygonal knots as its Hamiltonian cycles. In 1983 Conway and Gordon proved that any spatial embedding of $K_7$ contains at least one nontrivial knot as its Hamiltonian cycle \cite{CG}. And their result was generalized by Negami. He showed that for any knot type $\mathfrak{K}$, there exists a natural number $r(\mathfrak{K})$ such that every linearly embedded $K_{r(\mathfrak{K})}$ contains a polygonal knot of the type $\mathfrak{K}$ \cite{Negami}. For example $r(\mbox{trefoil})=7$ \cite{Alfonsin, Brown}. After Conway-Gordon's work, knots (also links) are considered as intrinsic properties of abstract graphs, and numerous subsequent works have been continued until recently. For a survey on this research  the readers are referred to \cite{Alfonsin2}.
In this paper we would like to study knots in $K_7$ more specifically.

\vspace{0.2cm}
Let $f:K_n \rightarrow \mathbb{R}^3$ be a linear embedding. From Theorem \ref{stick-thm} we know that any spatial trigon, tetragon and pentagon are trivial knots. Hence, for $n=3, 4$ and $5$, $f(K_n)$ contains only trivial knots as its cycles.

By Theorem \ref{stick-thm} any spatial hexagon is trivial or trefoil, and trefoil is the only nontrivial knot type which may appear in $f(K_6)$. Hence we may ask how many hexagonal trefoil knots can reside in $f(K_6)$ as Hamiltonial cycles. This question was answered in a previous work of the author \cite{HJ}, in which the point configuration of vertices of hexagonal trefoil knot was investigated and it was proved that {\em $f(K_6)$ contains at most one trefoil knot}. Recently, using the second coefficient of the Conway polynomial for knots, this was reproved by Nikkuni \cite{Nikkuni}.

\begin{figure}
{\epsfxsize=7cm \epsfbox{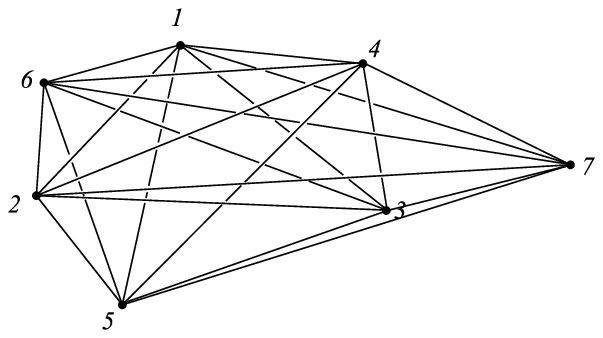}}
\caption{ }
\label{fig0}
\end{figure}

Again by Theorem \ref{stick-thm} any heptagon is trivial, trefoil or figure-8, and hence the figure-8 is the knot type which may be newly found in $f(K_7)$. And, considering polygonal index to be a measure for knot complexity, we may say that the figure-8 is the largest knot type which can reside in $f(K_7)$. In this paper we determine how many heptagonal figure-8 knots can be contained in $f(K_7)$. The following is the main result of this paper.
\begin{thm} \label{main-thm}
Any linear embedding of $K_7$ contains at most three heptagonal figure-8 knots as its Hamiltonian cycles.
\end{thm}
\noindent Note that a generic set of seven points in $\mathbb{R}^3$ determine a linear embedding of $K_7$. Hence the theorem can be restated as follows:
\begin{cor}
Any seven points in general position of $\mathbb{R}^3$ constitute at most three heptagonal figure-8 knots.
\end{cor}
\noindent Figure \ref{fig0} illustrates a realization of linearly embedded $K_7$ which contains three figure-8 knots $\langle1234567\rangle$, $\langle1236754\rangle$ and $\langle1276345\rangle$. Therefore the upper bound in Theorem \ref{main-thm} is strict.

The rest of paper will be devoted to the proof of Theorem \ref{main-thm}.

\section{Heptagonal figure-8 knot}
In this section we introduce a previous result of the author \cite{H} which would be a key lemma to prove Theorem \ref{main-thm}. Let $P$ be a heptagonal knot. We will call the line segments of $P$ {\em edges of $P$} and their end points {\em vertices}. The vertices of $P$ are assumed to be in general position, that is, any four vertices are not coplanar. And label the vertices of $P$ by $\{1, 2, \ldots , 7\}$ so that each vertex $i$ is connected to $i+1$ (mod $7$) by an edge of $P$, that is, a labeling of vertices is determined by a choice of base vertex and an orientation of $P$. Given a labelling of vertices let $\Delta_{i_1i_2i_3}$ denote the triangle formed by three vertices $\{i_1, i_2, i_3\}$, and $e_{jk}$ the line segment from the vertex $j$ to $k$. The relative position of such a triangle and a line segment will be represented via ``$\epsilon$'' which is defined in the below:
\begin{itemize}
\item[(i)] If $\Delta_{i_1i_2i_3} \cap e_{jk} = \emptyset$, then set $\epsilon(i_1i_2i_3,jk)=0$.
\item[(ii)] Otherwise,\\
$\epsilon(i_1i_2i_3,jk)= 1$ ({\em resp. $-1$}), when $(\overrightarrow{i_1i_2}\times \overrightarrow{i_2i_3})\cdot \overrightarrow{jk} > 0$ ({\em resp. $<0$}).
\end{itemize}

The tables in Theorem \ref{key-lemma} show the values of $\epsilon$ between triangles formed by three consecutive vertices and edges of $P$. If $\epsilon$ is zero, then the corresponding cell in the table is filled by ``$\times$''. Otherwise, we mark by ``$+$'' or ``$-$'' according to the sign of $\epsilon$. For example, if $P$ is of Type-I,  then $\epsilon(123,67)=0$ and \\
\centerline{ $(\epsilon(123,45),\epsilon(123,56),\epsilon(234,56))=(1,-1,-1) \quad \mbox{or} \quad (-1,1,1)$.}


\begin{thm} \cite{H} \label{key-lemma}
Let $P$ be a heptagonal knot such that its vertices are in general position. Then $P$ is figure-8 if and only if the vertices of $P$ can be labelled so that the polygon satisfies one among three  types Type-I, II and III. \\

\centering{ \small
\begin{tabular}{|c|c|c|c|}
\hline
 & 45 & 56 & 67 \\ \cline{2-4}
\raisebox{1.5ex}[0pt]{123} & $\pm$ & $\mp$ & $\times$ \\ \hline
 & 56 & 67 & 71 \\ \cline{2-4}
\raisebox{1.5ex}[0pt]{234} & $\mp$ & $\times$ & $\times$ \\ \hline
 & 67 & 71 & 12 \\ \cline{2-4}
\raisebox{1.5ex}[0pt]{345} & $\times$ & $\pm$ & $\times$ \\ \hline
 & 71 & 12 & 23 \\ \cline{2-4}
\raisebox{1.5ex}[0pt]{456} & $\pm$ & $\times$ & $\times$ \\ \hline
 & 12 & 23 & 34 \\ \cline{2-4}
\raisebox{1.5ex}[0pt]{567} & $\times$ & $\mp$ & $\times$ \\ \hline
 & 23 & 34 & 45 \\ \cline{2-4}
\raisebox{1.5ex}[0pt]{671} & $\mp$ & $\times$ & $\times$ \\ \hline
 & 34 & 45 & 56 \\ \cline{2-4}
\raisebox{1.5ex}[0pt]{712} & $\times$ & $\pm$ & $\times$ \\ \hline
\multicolumn{4}{c}{ } \\
\multicolumn{4}{c}{Type-I}
\end{tabular}
\hspace{0.5cm}
\begin{tabular}{|c|c|c|c|}
\hline
 & 45 & 56 & 67 \\ \cline{2-4}
\raisebox{1.5ex}[0pt]{123} & $\pm$ & $\mp$ & $\times$ \\ \hline
 & 56 & 67 & 71 \\ \cline{2-4}
\raisebox{1.5ex}[0pt]{234} & $\mp$ & $\times$ & $\times$ \\ \hline
 & 67 & 71 & 12 \\ \cline{2-4}
\raisebox{1.5ex}[0pt]{345} & $\times$ & $\pm$ & $\times$ \\ \hline
 & 71 & 12 & 23 \\ \cline{2-4}
\raisebox{1.5ex}[0pt]{456} & $\pm$ & $\times$ & $\times$ \\ \hline
 & 12 & 23 & 34 \\ \cline{2-4}
\raisebox{1.5ex}[0pt]{567} & $\times$ & $\mp$ & $\times$ \\ \hline
 & 23 & 34 & 45 \\ \cline{2-4}
\raisebox{1.5ex}[0pt]{671} & $\mp$ & $\pm$ & $\times$ \\ \hline
 & 34 & 45 & 56 \\ \cline{2-4}
\raisebox{1.5ex}[0pt]{712} & $\times$ & $\pm$ & $\times$ \\ \hline
\multicolumn{4}{c}{ } \\
\multicolumn{4}{c}{Type-II}
\end{tabular}
\hspace{0.5cm}
\begin{tabular}{|c|c|c|c|}
\hline
 & 45 & 56 & 67 \\ \cline{2-4}
\raisebox{1.5ex}[0pt]{123} & $\pm$ & $\mp$ & $\times$ \\ \hline
 & 56 & 67 & 71 \\ \cline{2-4}
\raisebox{1.5ex}[0pt]{234} & $\times$ & $\mp$ & $\times$ \\ \hline
 & 67 & 71 & 12 \\ \cline{2-4}
\raisebox{1.5ex}[0pt]{345} & $\times$ & $\pm$ & $\times$ \\ \hline
 & 71 & 12 & 23 \\ \cline{2-4}
\raisebox{1.5ex}[0pt]{456} & $\pm$ & $\times$ & $\times$ \\ \hline
 & 12 & 23 & 34 \\ \cline{2-4}
\raisebox{1.5ex}[0pt]{567} & $\times$ & $\mp$ & $\times$ \\ \hline
 & 23 & 34 & 45 \\ \cline{2-4}
\raisebox{1.5ex}[0pt]{671} & $\mp$ & $\times$ & $\times$ \\ \hline
 & 34 & 45 & 56 \\ \cline{2-4}
\raisebox{1.5ex}[0pt]{712} & $\times$ & $\pm$ & $\times$ \\ \hline
\multicolumn{4}{c}{ } \\
\multicolumn{4}{c}{Type-III}
\end{tabular}
}
\end{thm}

\section{Proof of Theorem \ref{main-thm}}
Let $K$ be a linearly embedded $K_7$ which contains a heptagonal figure-8 knot $P$ as its Hamiltonian cycle. Theorem \ref{key-lemma} guarantees that the vertices of $K$ can be labelled by $\{1,2,3,4,5,6,7\}$ so that $P$ corresponds to the cycle $\langle 1234567 \rangle$ and satisfies one among Type-I, II and III. Therefore our main theorem comes from the three lemmas in the below.
\begin{lem} \label{lemma-1}
If $P$ is of Type-I, then $K$ contains at most two heptagonal figure-8 knots as its Hamiltonian cycles.
\end{lem}
\begin{lem} \label{lemma-2}
If $P$ is of Type-II, then $K$ contains at most three heptagonal figure-8 knots as its Hamiltonian cycles.
\end{lem}
\begin{lem} \label{lemma-3}
If $P$ is of Type-III, then it is the only figure-8 knot among the Hamiltonian cycles of $K$.
\end{lem}

\begin{figure}
\centerline{\epsfbox{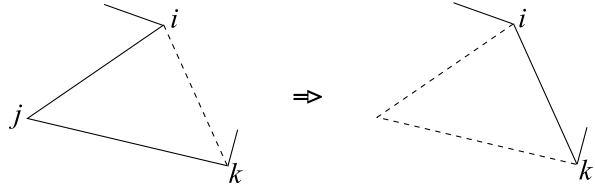}}
\caption{Reducing $\langle \cdots ijk \cdots \rangle$ into $\langle \cdots ik \cdots \rangle$ along $\Delta_{ijk}$.}
\label{fig2}
\end{figure}

These lemmas will be proved in next three sections. To prove them we utilize a concept, called {\em trivial triple}. A set of three vertices $\{i,j,k\}$ of $K$ will be called a {trivial triple}, if $\epsilon (ijk, lm)=0$ for every two vertices $l, m \notin \{i, j, k\}$. Suppose that a heptagon $Q$ in $K$ has the trivial triple as its consecutive three vertices. Then $Q$ is not figure-8, because we can isotope $Q$ to be a hexagon by reduction along $\Delta_{ijk}$ as illustrated in Figure \ref{fig2}.
In each case of Type-I, II and III we will try to find trivial triples as many as possible, so that heptagons which have any of such triples as consecutive vertices are excluded from possible Hamiltonian cycles of figure-8. And then the remaining Hamiltonian cycles will be observed more closely to check whether they are figure-8 or not.

\section{Proof of Lemma \ref{lemma-1}} Throughout this section we suppose that $P$ is of Type-I. Then, without loss of generality, it may be assumed that $\epsilon(123,45)=+1$. We begin the proof of Lemma \ref{lemma-1} with three claims in the below. Let $I_{ijk}$ denote the plane containing three vertices $\{i,j,k\}$. For an ordered sequence $ijk$ of the three vertices, define
$$ H^+_{ijk} = \{p\in \mathbb{R}^3 \; | \; (\overrightarrow{ij}\times \overrightarrow{jk})\cdot \overrightarrow{jp} > 0 \} \;\; \mbox{and} \;\;
H^-_{ijk} = \{q\in \mathbb{R}^3 \; | \; (\overrightarrow{ij}\times \overrightarrow{jk})\cdot \overrightarrow{jq} < 0 \} .$$

\begin{figure}
\centerline{\epsfbox{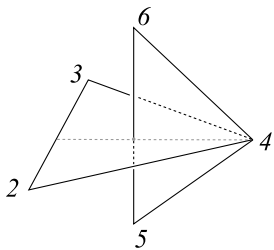}}
\caption{$\epsilon(234,56)=-1$ }
\label{fig3}
\end{figure}

\vspace{0.2cm}
\noindent
\textbf{Claim (i).} {\em The plane $I_{456}$ separates the two vertices $\{1,3\}$ from $\{2,7\}$.} \\
\noindent {\em Proof.} Under our assumption, $\epsilon(234,56)$ is $-1$, hence we know that $I_{456}$ separates the vertex $2$ from $3$ as illustrated in Figure \ref{fig3}. In fact the vertex $2$ ({\em resp. }3) belongs to $H^-_{456}$ ({\em resp.} $H^+_{456}$). Also from $\epsilon(456,71)=1$ we have that $7 \in H^-_{456}$ and $1 \in H^+_{456}$. $\Box$

\vspace{0.2cm}
From Claim (i) we can find out some trivial triples. For example the interior of $\Delta_{134}$ belongs to $H^+_{456}$, but the four vertices $\{2,5,6,7\}$ do not belong to the open half space. This implies that there is no edge of $K$ which penetrates $\Delta_{134}$. In this way we have a set of trivial triples in the below.
$$S_{1}= \{\{1, 3, 4\}, \{1, 3, 5\}, \{1, 3, 6\},\{2, 4, 7\}, \{2, 5, 7\}, \{2, 6, 7\}\}$$

\noindent
\textbf{Claim (ii).} {\em The plane $I_{671}$ separates the two vertices $\{3,5\}$  from $\{2,4\}$.} \\
\noindent
\textbf{Claim (iii).} {\em The plane $I_{123}$ separates the two vertices $\{4,6\}$ from $\{5,7\}$.} \\

Claim (ii) can be proved from $\epsilon(456,71)=1$ and $\epsilon(671,23)=-1$ in a similar way. Also a set of trivial triples is derived from the claim.
$$S_2=\{\{3, 5, 6\}, \{3, 5, 7\}, \{2, 4, 6\}, \{1, 2, 4\} \} $$

Claim (iii) is proved from $\epsilon(671,23)=-1$ and $\epsilon(123,45)=1$, and the third set of trivial triples is found out.
$$S_3=\{\{1, 5, 7\}, \{1, 4, 6\}, \{3, 4, 6\}\} $$

\begin{figure}
\centerline{\epsfbox{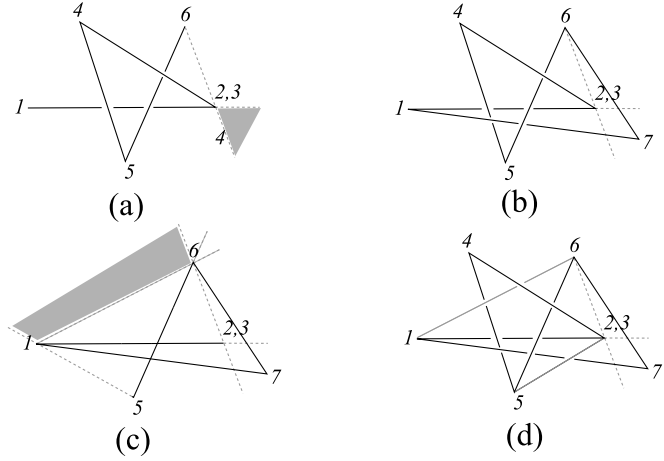}}
\caption{The image of $P$ under $\pi$}
\label{fig4}
\end{figure}

Let $N$ be a plane orthogonal to $\overrightarrow{23}$ and $\pi : \mathbb{R}^3 \equiv N\times\mathbb{R} \rightarrow N$ be the orthogonal projection onto $N$ such that the vertex $3$ is above the vertex $2$ with respect to the $\mathbb{R}$-coordinate. Then $\Delta_{123}$ is projected to a line segment. Now we observe the projected image of $P$ under $\pi$. Since $\epsilon(123,45)=1$ and $\epsilon(123,56)=-1$, the edges $e_{45}$ and $e_{56}$ should pass over $e_{12}$ with respect to the $\mathbb{R}$-coordinate. Also $\epsilon(234,56)=-1$ implies that $e_{34}$ passes over $e_{56}$. Figure \ref{fig4}-(a) shows the image of edges from $1$ to $6$ under $\pi$. Note that $\pi(\{7\})$ should be located in the shaded region, because the edge $e_{23}$ penetrates $\Delta_{671}$. Now we observe the relative position of $e_{71}$ with respect to $e_{45}$ and $e_{56}$. $\pi(\{1,7\})$ belongs to the outside of $\pi(\Delta_{456})$, and $e_{71}$ penetrates $\Delta_{456}$. This implies that $\pi(e_{71})$ should intersect both $\pi(e_{45})$ and $\pi(e_{56})$. On the other hand $e_{45}$ penetrates $\Delta_{712}$, and as shown in Figure \ref{fig4}-(a) the edge passes over $e_{12}$. Hence $e_{45}$ should pass under $e_{71}$. Consequently, for $e_{71}$ to penetrate $\Delta_{456}$, it should pass under $e_{56}$. From the resulting image of $\pi(P)$ in Figure \ref{fig4}-(b) we immediately find three more trivial triples
$$S_4=\{\{1, 4, 5\}, \{2, 3, 7\}, \{3, 6, 7\}\} \;\; .$$

Also from the image we know that $e_{47}$ is the only edge of $K$ which may penetrate $\Delta_{125}$. But this is not possible because $\epsilon(712,45)=1$, that is, $I_{712}$ separates the vertex $4$ from $5$. Hence $\{1, 2 ,5\}$ is another trivial triple.
From Claim (i) we know that $e_{17}$, $e_{37}$ and $e_{67}$ are the only edges of $K$ that may penetrate $\Delta_{245}$. But again the image of $P$ enables us to exclude $e_{37}$ and $e_{67}$ from the possibility. Also the edge $e_{17}$ should be excluded, because $e_{45}$ penetrates $\Delta_{712}$, that is, $\Delta_{245}\cap I_{712}$ belongs to $\Delta_{712}$. Therefore we have the fifth set of trivial triples
$$S_5=\{\{1, 2, 5\}, \{2,4,5\}\} \;\;.$$

To find more trivial triples we consider two subcases according to the position of the vertex $7$ with respect to the plane $I_{234}$. Figure \ref{fig4}-(b) shows that the vertices $\{1, 5\}$ belong to the half space $H^+_{234}$. If the vertex $7$ belongs to the other half space $H^-_{234}$, then $\{4,6,7\}$ is a trivial triple. Now suppose that the vertex belongs to $H^+_{234}$. Then $\{2,3,6\}$ is trivial, because the vertex $6$ is the only vertex in $H^-_{234}$. We will show that two more triples $\{1,4,7\}$ and $\{4,5,7\}$ are trivial in this case. From Claim (ii) we can know that $e_{23}$, $e_{25}$ and $e_{26}$ are the only edges of $K$ which may penetrate $\Delta_{147}$. But $e_{23}$ and $e_{26}$ are excluded from the candidates, because $\pi(\Delta_{147}) \cap \pi(\Delta_{236})$ is empty. Recall $\epsilon(712,45)=1$. Considering the relative positions of the five vertices we know that $e_{25}$ can not penetrate $\Delta_{147}$. Again, combining  $\pi(\Delta_{457}) \cap \pi(\Delta_{236})= \emptyset$ and $\epsilon(712,45)=1$ with Claim (i), we know that the triple $\{4,5,7\}$ is trivial. Therefore we have one more set of trivial triples, that is, either $S_6$ or $S_{6^{\prime}}$.
$$S_6=\{\{4, 6, 7\}\}, \;\;\; S_{6^{\prime}}= \{\{2, 3, 6\}, \{1, 4, 7\}, \{4,5,7\}\} \;\;$$

Now we observe the relative position of $e_{61}$ with respect to $e_{34}$ and $e_{45}$. A key information for this observation is $\epsilon(671,34)=0$ and $\epsilon(671,23)=-1$, that is, $e_{34}$ dose not penetrate $\Delta_{671}$ and the vertex $3$ belongs to $H^-_{671}$. Firstly note that $\pi(\{4\})$ should not be contained in $\pi(\Delta_{126})$. Otherwise, since $\pi(\Delta_{126}) \subset\pi(\Delta_{671})$, our key information implies that the vertices $4$ and $3$ belong to the same half space with respect to $I_{671}$. But it is contradictory to Claim (ii). Therefore $\pi(\{4\})$ should belong to the shaded region in Figure \ref{fig4}-(c), which implies that $\pi(e_{61})$ intersect both $\pi(e_{45})$ and $\pi(e_{34})$. In fact  $e_{34}$ should pass over $e_{61}$, because the vertex $3$ belongs to $H^-_{671}$, that is, $3$ is above $\Delta_{671}$ with respect the $\mathbb{R}$-coordinate, and $e_{34}$ does not penetrate the triangle. And  $e_{45}$ should pass under $e_{61}$, because the edge passes under $e_{71}$ but does not pierce $\Delta_{671}$. In addition it can be seen that $e_{71}$ passes under $e_{35}$, because the edge penetrates $\Delta_{345}$ and passes over $e_{45}$. Figure \ref{fig4}-(d) shows the resulting image $\pi(P\cup e_{61} \cup e_{35})$.

The observation in the above helps us find another trivial triple $\{1,5,6 \}$. From Claim (ii) it can be known that $e_{37}$, $e_{23}$ and $e_{34}$ are the only edges which may penetrate $\Delta_{156}$.
But looking into  Figure \ref{fig4}-(d) we can exclude the first two candidates. And also the third is excluded, because $e_{34}$ passes over both $e_{61}$ and $e_{56}$.

Let $S$ be the set of trivial triples of $P$ that we have found out so far. Then,
\begin{eqnarray*}
S & = & \big( \cup_{i=1}^{5} S_i \big) \cup S_{6^{\prime}} \cup \{1,5,6\} \;\;\;   \mbox{if} \; 7\in H^+_{234} \\
& & \big( \cup_{i=1}^{5} S_i \big) \cup S_6 \cup \{1,5,6\} \;\;\;   \mbox{if} \; 7\in H^-_{234}
\end{eqnarray*}

\begin{figure}
\centerline{\epsfbox{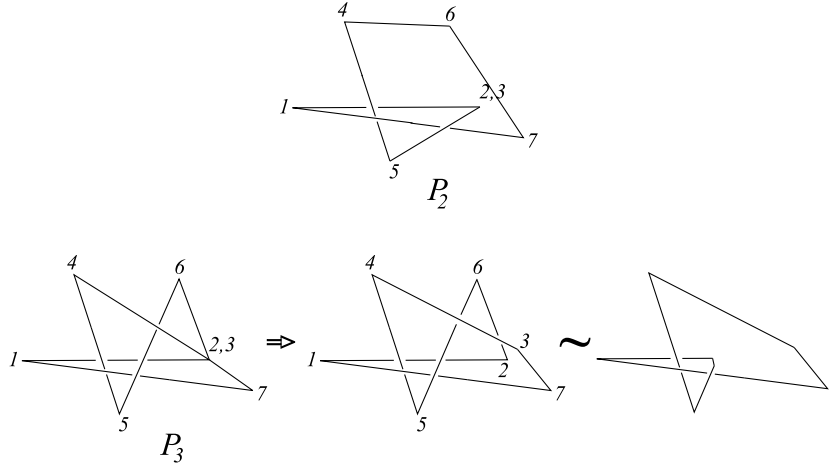}}
\caption{$P_2$ and $P_3$ }
\label{fig5}
\end{figure}

Now we exclude all Hamiltonian cycles of $K$ which contain any triple in $S$ as consecutive three vertices, because such cycles are not figure-8 as mentioned in the previous section. In the case that $7\in H^+_{234}$, there remain three heptagons:
$$ P_1=P=\langle1234567\rangle , \;\; P_2=\langle1235467\rangle, \;\; P_3=\langle1265437\rangle $$
Figure \ref{fig5} shows the images of $P_2$ and $P_3$ under $\pi$ which are constructed by using $\pi(P)$. As seen in the figure $P_2$ is a trefoil. To see that also $P_3$ is a trefoil, perturb the heptagon slightly. Then, since the vertex $3$ is above the vertex $2$, we obtain another image representing a trefoil knot.

\begin{figure}
\centerline{\epsfbox{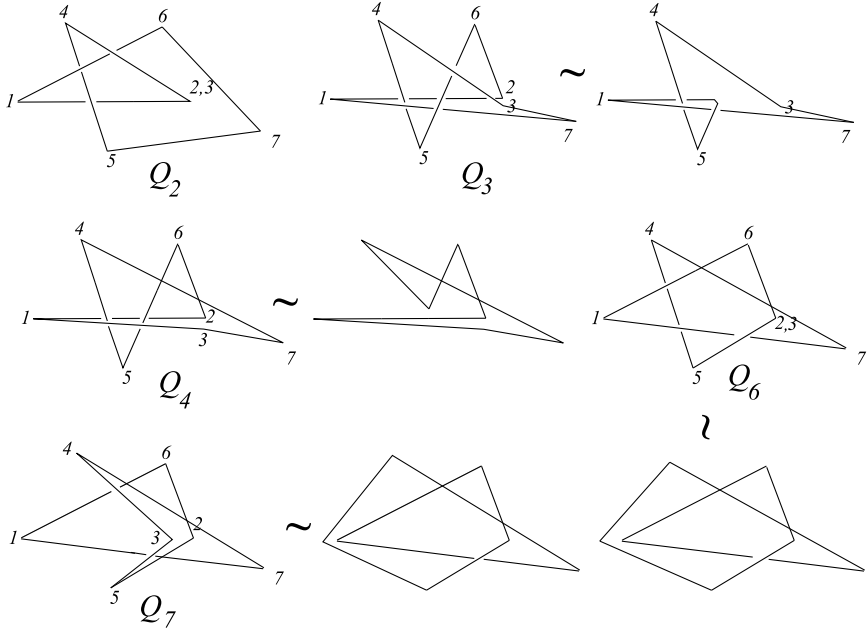}}
\caption{$Q_2$, $Q_3$, $Q_4$, $Q_6$ and $Q_7$ }
\label{fig6}
\end{figure}

In the case that $7\in H^-_{234}$, seven heptagons remain:
$$
\begin{array}{lll}
Q_1=P=\langle1234567\rangle & Q_2=\langle1234576\rangle & Q_3=\langle1265437\rangle\\
Q_4=\langle1265473\rangle & Q_5=\langle1326547\rangle & Q_6=\langle1623547\rangle\\
Q_7=\langle1625347\rangle &  &
\end{array}
$$
Figure \ref{fig6} shows their images under $\pi$. The polygons $Q_3$, $Q_4$ and $Q_7$ were perturbed slightly so that the vertices $2$ and $3$ are projected to different points. As shown in the figure $Q_2$ is a trefoil knot, and $Q_3$ can be isotoped to a polygon representing trefoil.
For $Q_4$ we need to recall $\epsilon(123,45)=1$ and $\epsilon(123,56)=-1$, which implies that $e_{31}$ passes over both $e_{45}$ and $e_{56}$. Hence $e_{45}$ and $e_{56}$ can be isotoped so that the resulting heptagon has exactly two double points with respect to $\pi$. Therefore $Q_4$ is a trivial knot. Also $Q_6$ can be isotoped to a octagon which has three double points, which implies that it is either a trefoil knot or a trivial knot. Similarly $Q_7$ is not a figure-8 knot.

The proof of Lemma \ref{lemma-1} is completed.

\section{Proof of Lemma \ref{lemma-2}}
In this section we suppose that $P$ is of Type-II. And $\epsilon(123,45)$ is assumed to be $+1$. Note that the table Type-II is almost identical with Type-I. The only difference is that  $\epsilon(671,34)$ is $1$ in Type-II. Therefore, following the procedure in the previous section, we can find out a set of trivial triples $T$ which is almost same with the set $S$. The only difference is that $T$ does not include $\{1,5,6\}$.
\begin{eqnarray*}
T & = & \big( \cup_{i=1}^{5} S_i \big) \cup S_{6^{\prime}} \;\;\;   \mbox{if} \; 7\in H^+_{234} \\
& & \big( \cup_{i=1}^{5} S_i \big) \cup S_6 \;\;\;\;   \mbox{if} \; 7\in H^-_{234}
\end{eqnarray*}

\begin{figure}
\centerline{\epsfbox{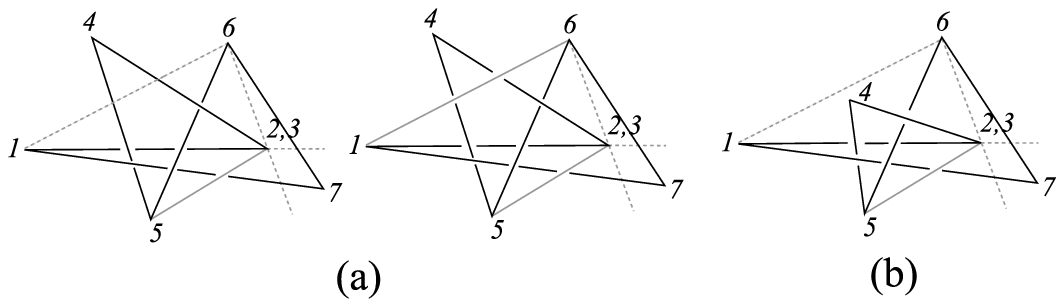}}
\caption{ }
\label{fig7}
\end{figure}

Also in the same way with the previous section, the image $\pi(P \cup e_{35})$ can be obtained. We observe the image more concretely. If $\pi(\{4\})$ does not belong to $\pi(\Delta_{126})$, then $\pi(e_{61})$ should intersect $\pi(e_{45})$ and $\pi(e_{34})$ as shown in the left of Figure \ref{fig7}-(a). Since $e_{45}$ does not penetrate $\Delta_{671}$ and passes under $e_{71}$, the edge should pass under $e_{61}$. In addition $e_{34}$ should pass under $e_{61}$ so that $\epsilon(671,34)$ is $1$, because $\epsilon(671,23)$ is $-1$, that is, the vertex $3$ is above the plane $I_{671}$ with respect to the $\mathbb{R}$-coordinate. The right of Figure \ref{fig7}-(a) shows the resulting image of $P \cup e_{35} \cup e_{61}$. Figure \ref{fig7}-(b) illustrates the case that $\pi(\{4\})$ belongs to $\pi(\Delta_{126})$.

In the case that the vertex $7$ belongs to $H^+_{234}$, after throwing away Hamiltonian cycles which have any triple in $T$ as consecutive three vertices, we have three remaining heptagons $\langle1234567\rangle$, $\langle1235467\rangle$ and $\langle1265437\rangle $.
Construct their images by using $\pi(P)$. Then, irrespective of the position of $\pi(\{4\})$, the resulting images are same with those in the previous section. Therefore the last two polygons can be shown to be trefoil knots.

\begin{figure}
\centerline{\epsfbox{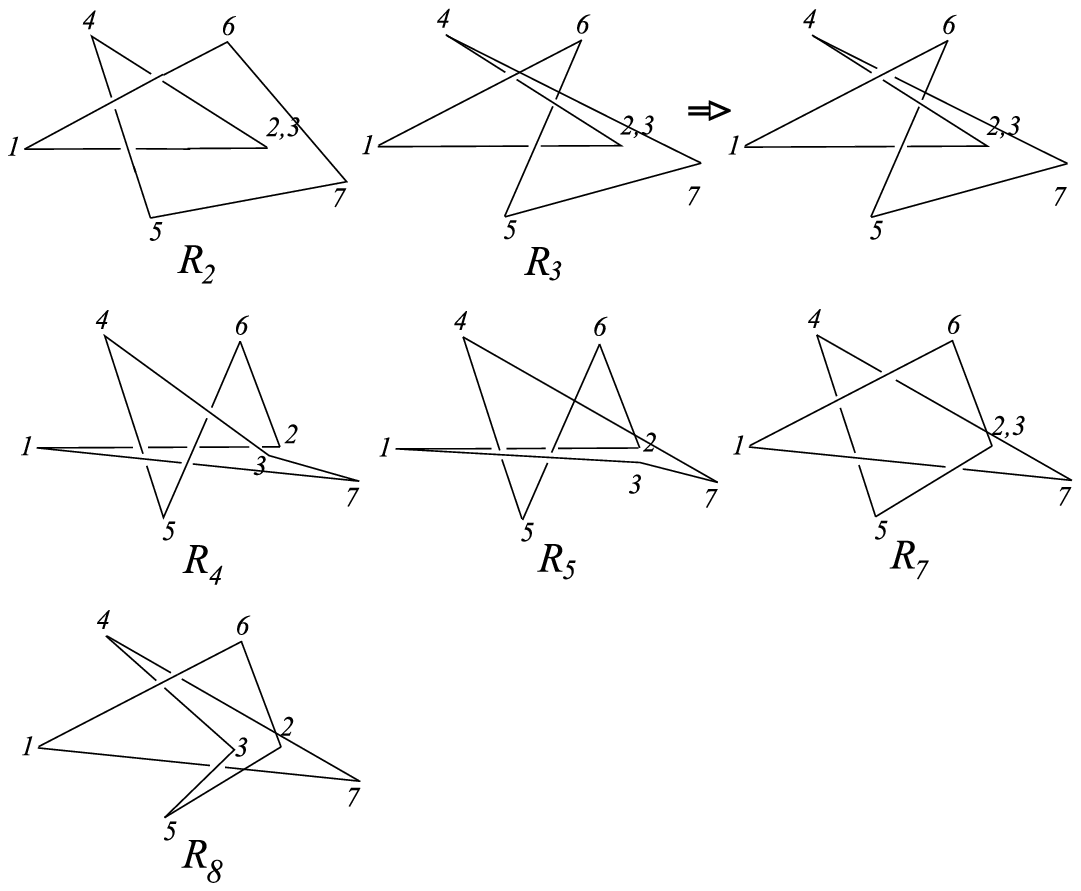}}
\caption{ }
\label{fig8}
\end{figure}

In the case that $H^-_{234}$ there remain nine heptagons:
$$
\begin{array}{lll}
R_1=P=\langle1234567\rangle & R_2=\langle1234576\rangle & R_3=\langle1234756\rangle\\
R_4=\langle1265437\rangle & R_5=\langle1265473\rangle & R_6=\langle1326547\rangle\\
R_7=\langle1623547\rangle & R_8=\langle1625347\rangle & R_9=\langle1652347\rangle
\end{array}
$$
Figure \ref{fig8} shows their images under $\pi$ which are constructed by using Figure \ref{fig7}-(a) after slight perturbation if necessary. It is easy to see that $R_2$, $R_4$ and $R_5$ are not figure-8. For $R_3$, $R_7$ and $R_8$ we need to recall that $\epsilon(671,34)$ is not zero. This implies that the edge $e_{61}$ does not penetrate $\Delta_{347}$. Therefore, since the edge passes over $e_{34}$, it should pass over $e_{47}$. From the resulting images we know that the three heptagons can be isotoped to other polygons which have exactly three double points. The images obtained from \ref{fig7}-(b) represent the same knot types with those from Figure \ref{fig7}-(a).

The proof of the lemma is completed.

\section{Proof of Lemma \ref{lemma-3}}
\begin{figure}
\centerline{\epsfbox{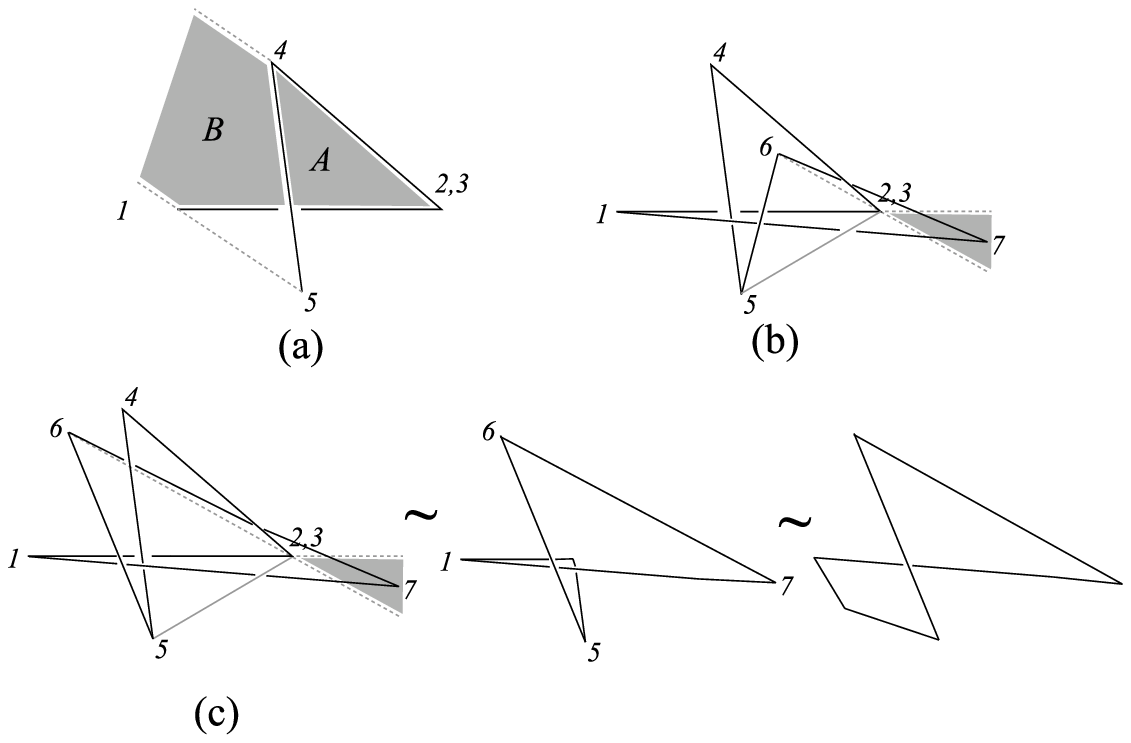}}
\caption{ }
\label{fig9}
\end{figure}

In this section $P$ is of Type-III and  $\epsilon(123,45)$ is assumed to be  $+1$. Using the table Type-III we attempt to construct the image of $P$ under the projection $\pi$. Since $\epsilon(123,45)=1$, the edges from $1$ to $5$ should be projected as illustrated in Figure \ref{fig9}-(a). Note that $\epsilon(123,56)=-1$ and $\epsilon(234,67)=-1$, that is, the vertex $6$ belongs to $H^-_{123} \cap H^+_{234}$. Hence the vertex $6$ should be projected into either the region $A$ or $B$.

Figure \ref{fig9}-(b) shows the image $\pi(P)$ when the vertex $6$ is projected into the region $A$. Since $\epsilon(671,23)$ is $-1$, the vertex $7$ should be projected into the shaded region. And since $\epsilon(456,71)$ is $1$, $e_{71}$ should intersect both $e_{45}$ and $e_{56}$. In fact the edge should pass over $e_{45}$, because $\epsilon(712,45)$ is $1$ and $e_{45}$ already passes over $e_{12}$. Combing this with $\epsilon(456,71)=1$ we know that $e_{71}$ passes under $e_{56}$. In addition the edge should pass under $e_{35}$, because $\epsilon(345,71)$ is $1$.

When the vertex $6$ projects into the region $B$, in a similar way, we can  construct the image $\pi(P)$ as illustrated in Figure \ref{fig9}-(c). One thing to be noted in this case is that $e_{67}$ should pass under $e_{45}$ because the edge passes under $e_{34}$ but $\epsilon(345,67)$ is $0$. Isotoping the resulting image it can be seen that $P$ is a trivial knot, a contradiction. Therefore this case can not happen.

From the image $\pi(P)$  we can find a set of trivial triples:
\begin{eqnarray*}
T_1 & = & \{ \{2,3,6\}, \{1,4,5\}, \{2,3,7\} \\
& & \{1,4,6\}, \{2,4,6\},\{3,4,6\}, \{1,5,7\},\{2,5,7\},\{3,5,7\}, \\
& & \{4,6,7\},\{2,4,7\},\{3,4,7\},\{1,5,6\},\{1,2,5\},\{1,3,5\}, \\
& & \{3,6,7\}, \{2,6,7\} \}
\end{eqnarray*}
See Figure \ref{fig9}-(b). It is obvious that $\{2,3,6\}$, $\{1,4,5\}$ and $\{2,3,7\}$ are trivial. Also the figure shows that $I_{123}$ separates $\{4,6\}$ from $\{5,7\}$, and $I_{236}$ separates $\{4,7\}$ from $\{1,5\}$, which implies that the next twelve triples are trivial. The edge $e_{24}$ is the only edge which may penetrate $\Delta_{367}$. For $\Delta_{267}$, $e_{34}$ is the only possible edge. But these possibilities are excluded, because $\epsilon(234,67)$ is $-1$.

\begin{figure}
\centerline{\epsfbox{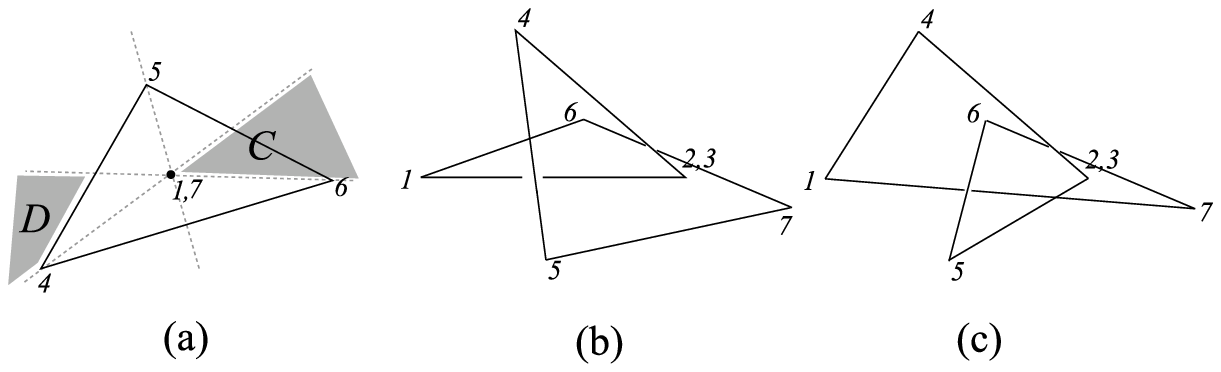}}
\caption{ }
\label{fig10}
\end{figure}

To find more trivial triples we observe another projected image of $P$. Let $N^{\prime}$ be a plane which is orthogonal to $\overrightarrow{17}$. And $\pi^{\prime}: N^{\prime}\times\mathbb{R} \rightarrow N^{\prime}$ is the projection map such that the vertex $7$ is above the vertex $1$ with respect to the $\mathbb{R}$-coordinate. Then $e_{71}$ is projected to a point in $\pi^{\prime}(\Delta_{456})$ as illustrated in Figure \ref{fig10}-(a). Since $\epsilon(671,23)$ is $-1$, the vertex $3$({\em resp.} $2$) is projected into  the upper({\em resp.} lower) half plane with respect to the line passing through $\pi^{\prime}(\Delta_{671})$. Hence $\epsilon(345,71)=1$ and $\epsilon(712,45)=1$ imply that $3$ belongs to the region $C$, and $2$ to $D$ respectively. From the figure four more trivial triples are found:
$$T_2=\{\{2,4,5\},\{1,3,6\},\{3,5,6\},\{1,2,4\} \}$$

After excluding all Hamiltonian cycles which have any of trivial triples in $T_1\cup T_2$ as consecutive three vertices, there remain only three heptagons:
$$U_1=P=\langle1234567\rangle, \;\; U_2=\langle1234576\rangle, \;\; U_3=\langle1432567\rangle $$
Figure \ref{fig10}-(b) and (c) show the images of $U_2$ and $U_3$ under $\pi$ respectively. It is clear that their knot types are not figure-8.

The proof is completed.

\end{document}